\documentclass[12pt,a4]{amsart}
\usepackage{amsmath}
\usepackage{amssymb}
\usepackage{amsthm}
\usepackage{enumitem}
\usepackage{url}
\usepackage{soul}
\usepackage{times}
\usepackage[T1]{fontenc}

\usepackage{color}
\usepackage{dsfont}
\usepackage{epsfig}
\usepackage[hidelinks]{hyperref}
\usepackage{setspace}
\usepackage{epstopdf}
\usepackage[authoryear,round,sort]{natbib}
\usepackage{comment}
\usepackage{booktabs}
\usepackage{float}
\usepackage{tikz}
\usepackage{algorithm}
\usepackage{stmaryrd}
\numberwithin{equation}{section}

\usepackage[font=small,labelfont=bf]{caption}
\DeclareCaptionFont{tiny}{\tiny}
\usepackage{geometry}
\theoremstyle{plain}

\newtheorem{theorem}{Theorem}[section]

\newtheorem{corollary}{Corollary}[section]

\theoremstyle{definition}

\newtheorem{example}{Example}[section]

\theoremstyle{remark}
\newtheorem{rk}{Remark}[section]
\expandafter\let\expandafter\oldproof\csname\string\proof\endcsname
\let\oldendproof\endproof
\renewenvironment{proof}[1][\proofname]{%
  \oldproof[\noindent\textbf{#1.} ]%
}{\oldendproof}
\newcommand{\1}{\mathds{1}}

\newcommand{\E}{\mathbb{E}}
\newcommand{\mP}{\mathbb{P}}
\newcommand{\be}{\begin{equation}}
\newcommand{\ee}{\end{equation}}
\newcommand{\by}{\begin{eqnarray*}}
\newcommand{\ey}{\end{eqnarray*}}

\renewcommand{\leq}{\leqslant}
\renewcommand{\geq}{\geqslant}
\usepackage{xcolor}
\definecolor{dark-red}{rgb}{0.4,0.15,0.15}
\definecolor{dark-blue}{rgb}{0.15,0.15,0.4}
\definecolor{medium-blue}{rgb}{0,0,0.5}
\hypersetup{
    colorlinks, linkcolor={red},
    citecolor={blue}, urlcolor={blue}
}

\begin{document}
\title{Hitting time and mixing time bounds of Stein's factors}
\author{Michael C.H. Choi}\thanks{}
\address{Department of Mathematics and Statistics, Hang Seng Management College, Hong Kong}
\email{michaelchoi@hsmc.edu.hk}
\date{\today}
\maketitle


\begin{abstract}
	For any discrete target distribution, we exploit the connection between Markov chains and Stein's method via the generator approach and express the solution of Stein's equation in terms of expected hitting time. This yields new upper bounds of Stein's factors in terms of the parameters of the Markov chain, such as mixing time and the gradient of expected hitting time. We compare the performance of these bounds with those in the literature, \textcolor{black}{and in particular we consider Stein's method for discrete uniform, binomial, geometric and hypergeometric distribution.}	As another application, the same methodology applies to bound expected hitting time via Stein's factors. This article highlights the interplay between Stein's method, modern Markov chain theory and classical fluctuation theory.
	\smallskip
	
	
\end{abstract}

\tableofcontents


\section{Introduction and main results}

Stein's method is well-known to be a powerful method for bounding the error rates of various distributional approximation, see e.g.  \cite{Ross11,LRS17,DH04,BC05} and the references therein. At the heart of it lies the Stein's equation
\begin{align}
	h(x) - \E h(Z) = Lf_h(x)\,,
\end{align}
where $Z \sim \pi$ is the target distribution with support on $\mathcal{X}$, $L$ is the Stein operator associated with $\pi$, $h$ belongs to a rich function class such as the class of indicator functions or Lipschitz continuous functions and $f_h$ is the solution of the Stein's equation. One popular approach to identify the Stein operator $L$ is the generator approach introduced by \cite{Barbour90,Gotze91}, where $L$ is the generator of a Markov process $X = (X_t)_{t \geq 0}$ with transition semigroup $(P_t)_{t \geq 0}$ on the state space $\mathcal{X}$ and stationary distribution $\pi$. Writing $\pi(f) = \int f d \pi$, the solution $f_h$ can then be related to $X$ via
\begin{align}\label{eq:steinsol}
f_h(x) = \int_0^{\infty} P_t h(x) - \pi(h)\, dt,
\end{align}
whenever the above integral exists. The obvious advantage of this approach is the connection with Markov processes, see e.g. \cite{BX01,ER08}. In addition, the solution form of $f_h$ naturally invites spectral techniques for Stein's method, which has been the subject of investigation in \cite{Schoutens01}. We also remark that \cite{DGV17} propose an interesting iterative procedure for bounding $f_h$ and its derivative.

Suppose now $\mathcal{X} = \llbracket 0,N \rrbracket$ with $N \in \mathbb{N}_0 \cup \{\infty\}$ is a countable set, where we write $\llbracket a,b \rrbracket = \{a,a+1,\ldots,b-1,b\}$ for $a,b \in \mathbb{Z}$. In Markov chain theory, the operator appearing on the right of \eqref{eq:steinsol}
$$D := \int_0^{\infty} P_t - \pi \, dt$$ is commonly known as deviation kernel $D = (D(i,j))_{i,j \in \mathcal{X}}$ as in \cite{Mao04,CvD02} (also known as fundamental matrix \cite{KSK}, ergodic potential \cite{Syski78} or centered resolvent \cite{Miclo16}). In this paper, we further exploit this intimate connection between Markov chain theory and Stein's method, thereby allowing us to express $f_h$ in terms of hitting time of an associated birth-death chain for any discrete target distribution $\pi$ as in \cite{ER08}, and from there connects Stein's method with modern Markov chain literature and offer universal bounds of Stein's factors in terms of quantities such as mixing time and eigentime. 

Before we discuss our main results, we fix our notations and revisit various parameters of countable Markov chains. We refer readers to \cite{AF14,LPW09,MT06} for in-depth account of these topics. For any probability measure $\mu,\nu$ with support on $\mathcal{X}$, the total variation distance between $\mu$ and $\nu$ is 
$$||\mu - \nu||_{TV} := \sup_{A \subset \mathcal{X}} |\mu(A) - \nu(A)| = \dfrac{1}{2} \sum_{j \in \mathcal{X}} |\mu(j) - \nu(j)|.$$
For $f$ on $\mathcal{X}$, we write $||f||_{\infty} := \sup_x |f(x)|$, the sup-norm of $f$. In this paper, we are primarily interested in the following parameters associated with an ergodic countable Markov chain $X = (X_t)_{t \geq 0}$:
\begin{itemize}
	\item Worst-case mixing time: for any $\epsilon > 0$, $$t_{mix}(\epsilon) := \inf\left\{t;~\sup_i ||P_t(i,\cdot) - \pi||_{TV} < \epsilon\right\}.$$
	\item Average hitting time and relaxation time: $$t_{av} := \sum_{i,j \in \mathcal{X}} \E_i(\tau_j)\pi(i)\pi(j) ,$$ where $\tau_A := \inf\left\{t;~X_t \in A\right\}$ is the first hitting time of $A$ by $X$, and as usual we write $\tau_j = \tau_{\{j\}}$. Note that for uniformly ergodic Markov chain, the eigentime identity \cite{AF14,Mao04,CuiMao10} is given by
	$$t_{av} = \sum_{i = 1}^{N} \dfrac{1}{\lambda_i} < \infty,$$
	where $\lambda_0 = 0 < \lambda_1 \leq \lambda_2 \leq \ldots$ are eigenvalues of $-L$. A closely related parameter is the relaxation time
	$$t_{rel} := 1/\lambda_1,$$
	and for finite reversible Markov chains we have $t_{mix}(\epsilon) \leq t_{rel} \log(1/\epsilon\pi_{min})$ with $\pi_{min} := \min_i \pi(i)$, see e.g. \cite[Theorem $20.6$]{LPW09}.
	\item Worst-case expected strong stationary time \cite{Aldous82}:
	$$t_{sst} := \sup_i \inf\{\E_i(T_i): \mP_i(X_{T_i} = j) = \pi(j) ~ \forall j \in \mathcal{X}\}.$$
	\item Worst-case expected hitting time of large set \cite{PS15,Oliveria12}: for $0 < \alpha < 1/2$,
	$$t_{hit}(\alpha) := \sup_{x,A \subset \mathcal{X}, \pi(A) \geq \alpha} \E_x(\tau_{A}).$$
	\item Worst-case expected deviation of hitting time to a single state \cite{Aldous82}:
	$$t_{dev} := \sup_{i,k} \sum_{j \in \mathcal{X}} |\E_i(\tau_j) - \E_k(\tau_j)|\pi(j).$$
\end{itemize}
Note that \cite{Oliveria12,PS15} show that for finite Markov chains $t_{mix}(1/4)$ and $t_{hit}(\alpha)$ are equivalent up to a constant depending on $\alpha$, while \cite{Aldous82} proves the equivalence (up to some universal constants) between $t_{mix}(1/4), t_{sst}, t_{dev}$ for reversible finite Markov chains. With the above notations in mind, we are now ready to state the main result of this paper:

\begin{theorem}[Hitting time as Stein's factors]\label{thm:main}
	Suppose that $\pi$ is a discrete target distribution on $\mathcal{X} = \llbracket 0,N \rrbracket$ with $N \in \mathbb{N}_0 \cup \{\infty\}$. The deviation kernel $D$ associated with $L$ and $\pi$ exists and is finite if and only if
	$$\sum_{i \in \mathcal{X}} \pi(i) \E_i(\tau_0) < \infty.$$
	In such case, if $L$ is reversible, $i,j \in \mathcal{X}$ and for any $h$ such that $\pi(|h|) < \infty$, 
	$$\sum_{i \in \mathcal{X}} \pi(i) |h|(i) \E_i(\tau_j) < \infty,$$
	we can write
	\begin{align}
		D(j,j) &= \pi(j) \sum_{i \in \mathcal{X}} \pi(i) \E_i(\tau_j), \label{eq:djj} \\
		D(i,j) &= D(j,j) - \pi(j) \E_i(\tau_j), \label{eq:dij} \\
		f_h(i) &= Dh(i) = \sum_{j \in \mathcal{X}} h(j) D(i,j), \label{eq:stein}\\
		\triangledown f_h(i) &:= f_h(i+1) - f_h(i) = \sum_{j \in \mathcal{X}} h(j) \pi(j) \left(\E_i(\tau_j) - \E_{i+1}(\tau_j)\right). \label{eq:gradstein}
	\end{align}
	In particular, if $h = \delta_j$, the Dirac mass at $j$, the sup-norm of the Stein factors are
	\begin{align*}
		||f_{\delta_j}||_{\infty} &\leq 2t_{av}, \\ 
		||\triangledown f_{\delta_j}||_{\infty} &= \pi(j) \sup_{i \in \mathcal{X}} |\E_i(\tau_j) - \E_{i+1}(\tau_j)| \leq t_{dev}.
	\end{align*}
	Note that we can always pick $L$ to be a birth-death process, and in such case the expected hitting time $\E_i(\tau_j)$ is readily computable and is expressed solely in terms of $\pi$, see Remark \ref{rk:bd}.
\end{theorem}

\begin{rk}\label{rk:bd}
	For any discrete distribution $\pi$ on $\mathcal{X}$, it is shown in \cite{ER08} that we can always pick $L$ to be the generator of a birth-death process with birth rate $b_i := L(i,i+1) = (i+1) \pi(i+1)/\pi(i)$ for $i \in \llbracket 0, N - 1\rrbracket$ and death rate $d_i := L(i,i-1) = i$ for $i \in \llbracket 1,N \rrbracket.$ According to \cite{Kijima97,CvD02,Mao04}, for $i,j \in \mathcal{X} = \llbracket 0,N\rrbracket$, 
	$$\E_i(\tau_j) = \sum_{k=i}^{j-1} \dfrac{1}{b_k\pi(k)} \pi(\llbracket 0,k \rrbracket) \1_{i < j} + \sum_{k=j}^{i-1} \dfrac{1}{b_k\pi(k)} \pi(\llbracket k+1,N \rrbracket) \1_{i \geq j},$$
	and so the discrete forward gradient $\triangledown f_{\delta_j}$ can be written as
	$$\triangledown f_{\delta_j}(i) = \pi(j) \left( \dfrac{1}{b_i\pi(i)} \pi(\llbracket 0,i\rrbracket)\1_{i < j} + \dfrac{1}{b_i\pi(i)} \pi(\llbracket i+1,N\rrbracket)\1_{i \geq j}\right).$$
	Another formula of $\E_i(\tau_j)$ involves differences in summation of eigenvalues, which is given by, for example for $i < j$,
	$$\E_i(\tau_j) = \sum_{k=0}^{j-1} \dfrac{1}{\lambda_k^{\llbracket 0,j-1\rrbracket}} - \sum_{k=0}^{i-1} \dfrac{1}{\lambda_k^{\llbracket 0,i-1\rrbracket}},$$
	where $(\lambda_k^{\llbracket 0,i-1\rrbracket})_{k = 0}^{i-1}$ are the non-zero eigenvalues of $-L$ restricted to $\llbracket 0,i-1 \rrbracket$, see e.g. \cite{Fill,GMZ12} and the references therein. In any case all these expressions are expressed in terms of a given target $\pi$.
\end{rk}

\begin{rk}
	It is tempting to think that $||f_{\delta_j}||_{\infty}$ equals to $ D(j,j)$. Yet, while $D(i,j) \leq D(j,j)$ for any $i,j$, it is unclear to the author whether $|D(i,j)|$ is less than or equal to $D(j,j)$.
\end{rk}

Theorem \ref{thm:main} reveals that hitting time and other Markov chain parameters as described are closely related to the structure and properties of $f_h$. In particular, $||\triangledown f_{\delta_j}||_{\infty} \leq t_{dev}$. This upper bound allows us to bound the Stein factors using various parameters: 
\begin{corollary}[Bounding Stein's factors via hitting and mixing time]\label{cor:main}
Suppose that $\pi$ is a discrete distribution with finite support on $\mathcal{X} = \llbracket 0,N \rrbracket$ and $N < \infty$, and $L$ is a reversible generator. Let $\mathcal{H} = \{h:\mathcal{X} \mapsto [0,1]\}$ and $i \in \mathcal{X}$, then
	\begin{align*}
\sup_{h \in \mathcal{H}}|f_h(i)| &\leq \sum_{j \in \mathcal{X}} |D(i,j)| \leq 2t_{av}, \\
\sup_{h \in \mathcal{H}} |\triangledown f_h(i)| &\leq \sum_{j \in \mathcal{X}} \pi(j) |\E_i(\tau_j) - \E_{i+1}(\tau_j)| \leq t_{dev} \leq   \begin{cases} 
10 t_{sst}, \\
5 t_{mix}(1/4), \\
5 t_{rel}\log(4/\pi_{min}), \\
5 C_{\alpha} t_{hit}(\alpha),
\end{cases}
\end{align*}
where $0 < \alpha < 1/2$ and $C_{\alpha} > 0$ is an universal constant depending only on $\alpha$.
\end{corollary}

\begin{rk}
	In practice, we argue that the relaxation time bound is perhaps the most useful of all as exact expressions or bounds on spectral gap $t_{rel} = 1/\lambda_1$ are readily available for many models. As a side note, we can also offer upper bounds involving log-Sobolev constant by bounding the mixing time with that, see \cite{DSC96}. This may perhaps yield tighter upper bound due to the double logarithm.
\end{rk}
The rest of this paper is organized as follows. In Section \ref{sec:proof}, we present the proofs of Theorem \ref{thm:main} and Corollary \ref{cor:main}. In Section \ref{sec:boundstein}, we illustrate our main results by detailing a few examples involving common distributions. As another application, we demonstrate a way to bound expected hitting time via Stein's factor in Section \ref{sec:boundhit}. 

\section{Proofs of the main results}\label{sec:proof}

\subsection{Proof of Theorem \ref{thm:main}}

\cite{CvD02,Mao04} show that $D$ exists and is finite if and only if 
$$\sum_{i \in \mathcal{X}} \pi(i) \E_i(\tau_j) < \infty.$$
for some $j$ and then for all $j$ by irreducibility. The expressions of $D(i,j)$ and $D(j,j)$ follow readily from \cite[equation $5.5$, $5.7$]{CvD02}. Define the $\alpha$-potential kernel $D^{\alpha} = (D^{\alpha}(i,j))_{i,j \in \mathcal{X}}$ by
$$D^{\alpha}(i,j) := \int_0^{\infty} e^{-\alpha t} \left(P_t(i,j) - \pi_j\right)\,dt.$$
Note that by dominated convergence theorem,
$$D^{\alpha}h(i) = \int_0^{\infty} e^{-\alpha t} \left(P_t h(i) - \pi(h)\right)\,dt.$$
Under the proposed assumptions on $h$ and reversibility of $L$, we use \cite[Lemma $5.1$]{CvD02} to arrive at 
$$f_h(i) = \int_0^{\infty} \left(P_t h(i) - \pi(h)\right)\,dt = \lim_{\alpha \rightarrow 0} D^{\alpha} h(i) = Dh(i) = \sum_{j \in \mathcal{X}} h(j)D(i,j).$$
To prove \eqref{eq:gradstein}, we see that
$$\triangledown f_h(i) = \sum_{j \in \mathcal{X}} h(j) \left(D(i+1,j) - D(i,j)\right) = \sum_{j \in \mathcal{X}} h(j) \pi(j) \left(\E_i(\tau_j) - \E_{i+1}(\tau_j)\right).$$
Now, we take $h = \delta_j$ and $f_{\delta_j}(i) = D(i,j)$. Using triangle inequality, we have
\begin{align*}
|f_{\delta_j}(i)| &= |D(i,j)| \leq D(j,j) + \pi(j) \E_i(\tau_j) \leq \sum_{j \in \mathcal{X}} D(j,j) + \sum_{j \in \mathcal{X}} \pi(j) \E_i(\tau_j) = t_{av} + t_{av} = 2t_{av},  
\end{align*}
where we use the fact that $\sum_j D(j,j) = t_{av}$ and the random target lemma \cite[Lemma $10.1$]{LPW09} for the second term. As for the gradient of $f_{\delta_j}$, it is straight forward from \eqref{eq:gradstein} that 
$$||\triangledown f_{\delta_j}||_{\infty} = \pi(j) \sup_{i \in \mathcal{X}} |\E_i(\tau_j) - \E_{i+1}(\tau_j)| \leq t_{dev}.
$$

\subsection{Proof of Corollary \ref{cor:main}}

To arrive at the first equation, we make use of \eqref{eq:stein}, $||h||_{\infty} \leq 1$ and triangle inequality to arrive at 
$$|f_h(i)| \leq \sum_{j \in \mathcal{X}} |D(i,j)| \leq \sum_{j \in \mathcal{X}} D(j,j) + \sum_{j \in \mathcal{X}} \pi(j) \E_i(\tau_j) = t_{av} + \sum_{j \in \mathcal{X}} \pi(j) \E_i(\tau_j).$$
Note that by the random target lemma \cite[Lemma $10.1$]{LPW09}, the second term is independent of $i$ which implies $\sum_{j \in \mathcal{X}} \pi(j) \E_i(\tau_j) = t_{av}$, and desired result follows. For the second set of equation, we apply \eqref{eq:gradstein} and $||h||_{\infty} \leq 1$ to yield
$$|\triangledown f_h(i)| \leq \sum_{j \in \mathcal{X}} \pi(j) |\E_i(\tau_j) - \E_{i+1}(\tau_j)| \leq \sup_k \sum_{j \in \mathcal{X}} \pi(j) |\E_i(\tau_j) - \E_{k}(\tau_j)| \leq t_{dev}.$$
$t_{dev} \leq 10 t_{sst}$ follows from \cite[Lemma $16$]{Aldous82}, and making use of that together with \cite[Lemma $12$]{Aldous82} give
$$t_{dev} \leq 10 t_{sst} \leq 5t_{mix}(1/4).$$
Utilizing the above and \cite[Theorem $20.6$]{LPW09} leads to 
$$t_{dev} \leq 5t_{mix}(1/4) \leq 5 t_{rel} \log(4/\pi_{min}).$$
Finally, the main result in \cite{Oliveria12,PS15} gives 
$$t_{dev} \leq 5t_{mix}(1/4) \leq 5C_{\alpha}t_{hit}(\alpha).$$
\section{Bounding Stein's factors via hitting and mixing time - examples}\label{sec:boundstein}

In this section, we will discuss in detail several examples to illustrate both Theorem \ref{thm:main} and Corollary \ref{cor:main}, and compare with existing bounds in the literature. Our primary comparison is the bounds in \cite{ER08}. Notable results of this section are the $O(\log n)$ bound in Example \ref{ex:bin}, and bounds in Example \ref{ex:hypergeom}.

\begin{example}[discrete uniform on $\llbracket 0,n-1 \rrbracket$ with $n < \infty$]
In our first example, we look at $\mathcal{X}= \llbracket 0,n-1 \rrbracket$ and $\pi(i) = 1/n$ for $i \in \mathcal{X}$, and we take $L$ to be the uniform chain with $L(i,j) = 1$ for all $i \neq j$ as in \cite[Example $48$]{Aldous82}, which is reversible. The nice feature about this chain is that
$$\E_i(\tau_j) =  
\begin{cases}
	n-1 & \text{if } i \neq j, \\
	0 & \text{if } i = j.
\end{cases} $$
Also, note that the eigenvalues of $-L$ are $0$ with multiplicity $1$ and $1/n$ with multiplicities $n-1$. As a result, Theorem \ref{thm:main} now reads, for any $i,j \in \mathcal{X}$,
\begin{align*}
	D(j,j) &= \left(\dfrac{n-1}{n}\right)^2, \\
	D(i,j) &= - \dfrac{n-1}{n^2}, \\
	t_{av} &= \dfrac{n-1}{n}, \\
	||f_{\delta_j}||_{\infty} &= \left(\dfrac{n-1}{n}\right)^2 \leq 2 t_{av}, \\
	||\triangledown f_{\delta_j}||_{\infty} &= \dfrac{n-1}{n}.
\end{align*}
As for Corollary \ref{cor:main}, since $t_{av} = O(1/n)$, $$\sup_{h \in \mathcal{H}}||f_h||_{\infty} = O(1/n).$$ In addition, note that the relaxation time bound yields a crude upper bound of size $O(n \log n)$ for the gradient of Stein's solution:
$$\sup_{h \in \mathcal{H}} |\triangledown f_h(i)| \leq \dfrac{2(n-1)}{n} \leq \dfrac{5}{\lambda_1}\log(4/\pi_{min}) = 5n\log(4n).$$ 
\end{example}
\textcolor{black}{Note that Stein's method for discrete uniform distribution has first been considered in \citet[Chapter $2$]{DH04}. In particular, in Theorem $2.2.1$ of \citet{DH04}, it is shown in the proof that, in our notations, for odd $n$ and $h = \delta_S$, $||f_h||_{\infty} \leq (n-1)/2$, so our uniform bound of size $O(1/n)$ seems to be tighter in this setting.
}
\begin{example}[Binomial distribution on $\llbracket 0,n \rrbracket$ with parameters $n$ and $0 < p < 1$]\label{ex:bin}
In the second example, we consider $\pi(i) = C^n_i p^i(1-p)^{n-i}$ for $i \in \mathcal{X} = \llbracket 0,n \rrbracket$. \textcolor{black}{Stein's method for binomial distribution has also been considered in \citet{Ehm91}.} We take $L$ to be a birth-death process with birth rate $b_i = p(n-i)$ and death rate $d_i = (1-p)i$. The non-zero eigenvalues of $-L$ are $\lambda_k = k$ for $k \in \llbracket 1,n \rrbracket$ (see e.g. \cite{Schoutens00}), and it follows from Corollary \ref{cor:main} that 
\begin{align*}
t_{av} &= \sum_{k=1}^n \dfrac{1}{k} = O(\log n), \\
\sup_{h \in \mathcal{H}}||f_h||_{\infty} &= O(\log n).
\end{align*}
As for the relaxation time bound, we have
\begin{align*}
\sup_{h \in \mathcal{H}} ||\triangledown f_h||_{\infty} \leq 5 n\min \left\{\log \left(\dfrac{4}{1-p}\right),\log \left(\dfrac{4}{p}\right)\right\}.
\end{align*}
This $O(n)$ bound does not seem to be useful at all when compared with the $O(1)$ uniform bound of $\sup_{h \in \mathcal{H}} ||\triangledown f_h||_{\infty} \leq \min\{1/(1-p),1/p\}$ as in \cite{ER08}.
\end{example}

\begin{example}[Hypergeometric distribution on $\llbracket 0,r \rrbracket$ with parameters $n$, $r$ and $0< 2r \leq n$]\label{ex:hypergeom}
	In this example, we study the hypergeometric distribution 
	$$\pi(i) = \dfrac{C^r_i C^{n-r}_{r-i}}{C^n_r}$$
	for $i \in \llbracket 0,r \rrbracket$, and pick $L$ to be the generator of the Bernoulli-Laplace model, that is, it is a birth-death chain with birth rate $b_i$ and death rate $d_i$ to be respectively
	$$b_i = \dfrac{(r-i)(n-r-i)}{r(n-r)},~d_i = \dfrac{i^2}{r(n-r)}.$$
	The eigenvalues of $-L$ are 
	$$\lambda_i = \dfrac{i(n-i+1)}{r(n-r)},$$
	so Corollary \ref{cor:main} now reads
		\begin{align*}
	\sup_{h \in \mathcal{H}}|f_h(i)| &\leq 2r(n-r) \sum_{i=1}^r \dfrac{1}{i(n-i+1)} = \dfrac{2r(n-r)}{n}\left(\log r + O(1)\right), \\
	\sup_{h \in \mathcal{H}} |\triangledown f_h(i)| &\leq 
	\dfrac{5r(n-r)}{n} \log\left(\max_{i \in \llbracket 0,r\rrbracket}\dfrac{4 C^n_r}{C^r_i C^{n-r}_{r-i}}\right), 
	\end{align*}
	where the equality follows from \cite[Page $2114$]{DSC06}. \textcolor{black}{Existing work on Stein's method for hypergeometric distribution include \citet{RS98}, \citet[Section $4$]{R05} and \citet[Section $4$ Example $5$]{Schoutens00}, however in these work bounds for the Stein's factors cannot be found. We also adapt a different Stein's equation, namely the generator of the Bernoulli-Laplace model, than the one in existing literature.}
\end{example}
	
\begin{example}[Geometric distribution with success probability $0 < p < 1$]
	In this example, we use only estimates and information on hitting time to bound the Stein's factors for geometric distribution. More specifically, we look at $\pi(i) = (1-p)^ip$ for $i \in \mathcal{X} = \mathbb{N}_0$, and choose $L$ with unit per capita death rate and birth rate $b_i = (i+1)(1-p)$. It follows from Remark \ref{rk:bd} that for any fixed $j \in \mathcal{X}$, if $i < j$ then
	$$\E_i(\tau_j) \leq \E_0(\tau_j) = \sum_{k=0}^{j-1} \dfrac{1-(1-p)^{k+1}}{(k+1)(1-p)^{k+1}p} \leq \sum_{k=0}^{j-1} \dfrac{1}{(1-p)^{k+1}p}$$
	For any $i,j \in \mathcal{X}$, and using  Remark \ref{rk:bd}, we have
	\begin{align*}
		|\E_i(\tau_j) - \E_{i+1}(\tau_j)| \leq \dfrac{1}{b_i \pi(i)} \pi(\llbracket 0,i \rrbracket) \1_{i < j} + \dfrac{1}{b_i \pi(i)} \pi(\llbracket i+1,\infty \rrbracket) \1_{i \geq j}.
	\end{align*}
	Specializing into the case of geometric distribution leads to
	\begin{align*}
		\pi(j)|\E_i(\tau_j) - \E_{i+1}(\tau_j)| \leq (1-p)^{j-i-1} \dfrac{1-(1-p)^{i+1}}{i+1} \1_{i < j} + \dfrac{(1-p)^{j}}{i+1} \1_{i \geq j}.
	\end{align*}
	Now, using Corollary \ref{cor:main} and summing over all possible $j$ gives
	$$\sup_{h \in \mathcal{H}} |\triangledown f_h(i)| \leq \dfrac{1-(1-p)^{i+1}}{p(i+1)} + \dfrac{1-(1-p)^{i+1}}{p(i+1)} = \dfrac{2(1-(1-p)^{i+1})}{p(i+1)} \leq \dfrac{2}{p},$$	
	and so $\sup_{h \in \mathcal{H}} ||\triangledown f_h||_{\infty} \leq 2/p$. Note that in \cite{ER08} they obtain
	\begin{align*}
		\sup_{h \in \mathcal{H}} |\triangledown f_h(i)| &\leq \min\left\{\dfrac{1}{i},\dfrac{1+p}{i+1}\right\}, \\
		\sup_{h \in \mathcal{H}} ||\triangledown f_h||_{\infty} &\leq \min\{1,1+p\} \leq \dfrac{2}{p},
	\end{align*}
	so our bound is looser than existing bound.
\end{example}

\section{Bounding expected hitting time via Stein's factors}\label{sec:boundhit}

In the previous section, we have illustrated how we can use information of hitting time and mixing time to give bounds on Stein's factors. We aim at achieving the opposite in this section and illustrate how we can obtain estimate on the expected hitting time to $0$ of a Galton-Watson with immigration (GWI) process. 

Recall that the generator $L$ of GWI is a birth-death chain with $b_i = p(r+i)$, $d_i = i$ and $\pi$ being the negative binomial distribution with parameters $0 < p < 1$ and $r >0$. Taking $h = \delta_0$ in Theorem \ref{thm:main}, for $i \in \mathbb{N}_0$ we have
$$\triangledown f_{\delta_0}(i) = \pi(0) \left(\E_i(\tau_0) - \E_{i+1}(\tau_0)\right) = (1-p)^r \left(\E_i(\tau_0) - \E_{i+1}(\tau_0)\right) .$$
Next, we bound the above expression using the Stein's factor bound in \cite{BGX15} to arrive at
\begin{align*}
(1-p)^r \left|\E_i(\tau_0) - \E_{i+1}(\tau_0)\right| \leq \dfrac{1}{1-p} =: C,
\end{align*}
and the above yields the following linear bound:
$$\E_i(\tau_0) \leq C(1-p)^{-r}i.$$
This bound can perhaps be refined by adapting non-uniform bounds of Stein's factor.

\section*{Acknowledgements}
The author would like to thank the anonymous referee for constructive comments that improve the presentation of the paper.

\bibliographystyle{abbrvnat}
\bibliography{thesis}

\end{document}